\documentclass [reqno]{amsart}
\usepackage {latexsym,epsfig}
\usepackage [latin1]{inputenc}

\newcommand {\OO}{{\mathcal O}}
\newcommand {\CC}{{\mathbb C}}
\newcommand {\NN}{{\mathbb N}}
\newcommand {\PP}{{\mathbb P}}
\newcommand {\QQ}{{\mathbb Q}}

\newcommand {\dt}{\partial_t}
\newcommand {\dz}{\partial_z}
\newcommand {\di}{\partial_1}
\renewcommand {\dj}{\partial_2}
\newcommand {\dci}{\partial_{c_i}}
\newcommand {\dxi}{\partial_{x_i}}
\newcommand {\at}{|_{c_i=0}}

\DeclareInputMath {183}{\cdot}

\newcommand {\vdim}{\mathop {\rm vdim}\,}
\newcommand {\codim}{\mathop {\rm codim}\,}
\newcommand {\preprint}[2]{preprint \discretionary {#1/}{#2}{#1/#2}}

\newtheorem {theorem}{Theorem}[section]
\newtheorem {lemma}[theorem]{Lemma}

\newtheorem {corollary}[theorem]{Corollary}
\theoremstyle {definition}
\newtheorem {definition}[theorem]{Definition}
\newtheorem {example}[theorem]{Example}
\theoremstyle {remark}
\newtheorem {remark}[theorem]{Remark}

\begin {document}

  % Main file: latex main

\title [Relative Gromov-Witten invariants and the mirror formula]{%
  Relative Gromov-Witten invariants and the mirror formula}
\author {Andreas Gathmann}
\address {Harvard University, Department of Mathematics, Science Center,
          1 Oxford Street, Cambridge, MA 02138, USA}
\email {andreas@math.harvard.edu}
\thanks {Funded by the DFG scholarships Ga 636/1--1 and Ga 636/1--2.}

\begin {abstract}
  Let $X$ be a smooth complex projective variety, and let $ Y \subset X $ be a
  smooth very ample hypersurface such that $ -K_Y $ is nef. Using the technique
  of relative Gromov-Witten invariants, we give a new short and geometric proof
  of (a version of) the ``mirror formula'', i.e.\ we show that the generating
  function of the genus zero 1-point Gromov-Witten invariants of $Y$ can be
  obtained from that of $X$ by a certain change of variables (the so-called
  ``mirror transformation''). Moreover, we use the same techniques to give a
  similar expression for the (virtual) numbers of degree-$d$ plane rational
  curves meeting a smooth cubic at one point with multiplicity $ 3d $, which
  play a role in local mirror symmetry.
\end {abstract}

\maketitle

  % Main file: latex main

For a smooth very ample hypersurface $Y$ of a smooth complex projective variety
$X$, the theory of relative Gromov-Witten invariants gives rise to an algorithm
that allows one to compute the genus zero Gromov-Witten invariants of $Y$ from
those of $X$ \cite {Ga}. The goal of this paper is to show that in the case
when $ -K_Y $ is nef, this algorithm can be ``solved'' explicitly to obtain a
formula that expresses the generating function of the 1-point Gromov-Witten
invariants of $Y$ in terms of that of $X$. This so-called ``mirror formula''
(also denoted ``quantum Lefschetz hyperplane theorem'' by some authors) has
already been known for some time (\cite {Gi}, \cite {LLY}, \cite {K}, \cite
{B}, \cite {L}). Our approach however is entirely different and essentially
``elementary'' in the sense that it does not use any of the special techniques
that have been used in the previous proofs, like e.g.\ torus actions,
equivariant cohomology, or moduli spaces other than the usual spaces of stable
maps to $X$ and their subspaces. This does not only make our proof much simpler
than the previous ones, but also hopefully easier to generalize, e.g.\ to more
general hypersurfaces, or to higher genus of the curves.

Let us briefly recall the ideas and results from \cite {Ga}. For $ n \ge 0 $
and a homology class $ \beta \in H_2 (X) / $torsion we denote by $ \bar M_n
(X,\beta) $ the moduli space of $n$-pointed genus zero stable maps to $X$ of
class $ \beta $. For any $ m \ge 0 $ there are closed subspaces $ \bar M_{(m)}
(X,\beta) $ of $ \bar M_1 (X,\beta) $ that can be thought of as parametrizing
1-pointed rational curves in $X$ having multiplicity (at least) $m$ to $Y$ at
the marked point. (For simplicity, we suppress in the notation the dependence
of these spaces on $Y$.) These moduli spaces have expected codimension $m$ in $
\bar M_1 (X,\beta) $. In fact, they come equipped with natural virtual
fundamental classes $ [ \bar M_{(m)} (X,\beta) ]^{virt} $ of this expected
dimension. If $X$ is a projective space and $Y$ a hyperplane, then these moduli
spaces do have the expected dimension, and their virtual fundamental classes
are equal to the usual ones.

The idea is now to raise the multiplicity $m$ of the curves from $0$ up to $
Y·\beta+1 $ by one at a time. Curves with multiplicity (at least) $0$ are just
unrestricted curves in $X$, whereas a multiplicity of $ Y·\beta+1 $ forces at
least the irreducible curves to lie inside $Y$. In other words, we consider the
chain of inclusions
  \[ \bar M_1 (Y,\beta) \subset \bar M_{(Y·\beta)} (X,\beta) \subset
       \bar M_{(Y·\beta-1)} (X,\beta) \subset \cdots \subset
       \bar M_{(0)} (X,\beta) = \bar M_1 (X,\beta) \]
of ``virtual codimension one''. The main theorem of \cite {Ga} describes each
of these inclusions explicitly in terms of intersection theory. This gives us a
way to describe $ \bar M_1 (Y,\beta) $ inside $ \bar M_1 (X,\beta) $, and hence
to compute Gromov-Witten invariants of $Y$ in terms of those of $X$.

It is easy to write down a naïve guess what these inclusions should look like.
A stable map in $X$ has multiplicity at least $m$ to $Y$ if and only if the $
(m-1) $-jet of $ ev^* Y $ vanishes, where $ ev: \bar M_1 (X,\beta) \to X $
denotes the evaluation map. Hence the cycle $ \bar M_{(m+1)} (X,\beta) $ inside
$ \bar M_{(m)} (X,\beta) $ should just be the first Chern class of the line
bundle of $m$-jets modulo $ (m-1) $-jets of $ ev^* \OO (Y) $. This Chern class
is easily computed to be $ ev^* Y + m\psi $, where $ \psi $ is the ``cotangent
line class'', i.e.\ the first Chern class of the line bundle whose fiber at a
stable map $ (C,x,f) $ is the cotangent space of $C$ at the point $x$.

However, our above informal description of $ \bar M_{(m)} (X,\beta) $ as the
space of curves with multiplicity at least $m$ to $Y$ at the marked point
breaks down at the ``boundary'', i.e.\ at those curves where the marked point
lies on a component of the curve that lies completely inside $Y$, so that the
multiplicity becomes ``infinite''. Hence the above calculation receives
correction terms from these curves. Their explicit form is given by the
following theorem (see \cite {Ga} theorem 2.6).

\begin {theorem} \label {old-thm}
  For all $ m \ge 0 $ we have
    \[ (ev^* Y + m \psi) · [\bar M_{(m)} (X,\beta)]^{virt}
      = [\bar M_{(m+1)} (X,\beta)]^{virt}
      + [D_{(m)} (X,\beta)]^{virt}. \]
  Here, the correction term $ D_{(m)} (X,\beta) = \coprod_r \coprod_{B,M}
  D(X,B,M) $ is a disjoint union of individual terms
    \[ D(X,B,M) := \bar M_{1+r} (Y,\beta^{(0)}) \times_{Y^r}
         \prod_{i=1}^r \bar M_{(m^{(i)})} (X,\beta^{(i)}) \]
  where $ r \ge 0 $, $ B = (\beta^{(0)},\dots,\beta^{(r)}) $ with $ \beta^{(i)}
  \in H_2(X) / $torsion and $ \beta^{(i)} \neq 0 $ for $ i>0 $, and $ M =
  (m^{(1)},\dots,m^{(r)}) $ with $ m^{(i)} > 0 $. The maps to $ Y^r $ are the
  evaluation maps for the last $r$ marked points of $ \bar M_{1+r}
  (Y,\beta^{(0)}) $ and each of the marked points of $ \bar M_{(m^{(i)})}
  (X,\beta^{(i)}) $, respectively. The union in $ D_{(m)} (X,\beta) $ is taken
  over all $r$, $B$, and $M$ subject to the following three conditions:
  \begin {align*}
    \sum_{i=0}^r \beta^{(i)} = \beta
      \qquad &\mbox {(degree condition),} \\
    Y·\beta^{(0)} + \sum_{i=1}^r m^{(i)} = m
      \qquad &\mbox {(multiplicity condition),} \\
    \mbox {if $ \beta^{(0)} = 0 $ then $ r \ge 2 $}
      \qquad &\mbox {(stability condition).}
  \end {align*}
  In the equation of the theorem, the virtual fundamental class of the summands
  $ D(X,B,M) $ is defined to be $ \frac {m^{(1)} \cdots m^{(r)}}{r!} $ times
  the class induced by the virtual fundamental classes of the factors $ \bar
  M_{1+r} (Y,\beta^{(0)}) $ and $ \bar M_{(m^{(i)})} (X,\beta^{(i)}) $. The
  spaces $ D(X,B,M) $ can be considered to be subspaces of $ \bar M_1 (X,\beta)
  $ (see below), so the equation of the theorem makes sense in the Chow group
  of $ \bar M_1 (X,\beta) $.
\end {theorem}

Geometrically speaking, the moduli spaces $ D(X,B,M) $ in the correction terms
describe curves with $ r+1 $ irreducible components $ C^{(0)},\dots,C^{(r)} $
with homology classes $ \beta^{(0)},\dots,\beta^{(r)} $, such that $ C^{(0)} $
lies inside $Y$, and the $ C^{(i)} $ for $ i>0 $ intersect $ C^{(0)} $ in a
point where they have multiplicity $ m^{(i)} $ to $Y$. The marked point is
always on the component $ C^{(0)} $. Using this description, the spaces $
D(X,B,M) $ can be considered as subspaces of $ \bar M_1 (X,\beta) $. The
multiplicity condition ensures that they are actually subspaces of $ \bar
M_{(m)} (X,\beta) $ and have the correct expected dimension. The factor $ \frac
1 {r!} $ in the definition of the virtual fundamental class of the correction
terms is just combinatorial and corresponds to the choice of order of the
components $ C^{(1)},\dots,C^{(r)} $. In contrast, the factor $ m^{(1)} \cdots
m^{(r)} $ is of geometric nature and somewhat tricky to derive.

As an example of the theorem, consider the case where $ X=\PP^3 $, $ Y=H $ is a
hyperplane, and $ \beta $ is the class of cubic curves in $X$. Then the
equations of the theorem for $ m=0,\dots,3 $ can be pictured as follows (where
we set $ \bar M_{(m)} := \bar M_{(m)} (\PP^3,3) $):

\vspace {0mm plus2mm}
\begin {center} \epsfig {file=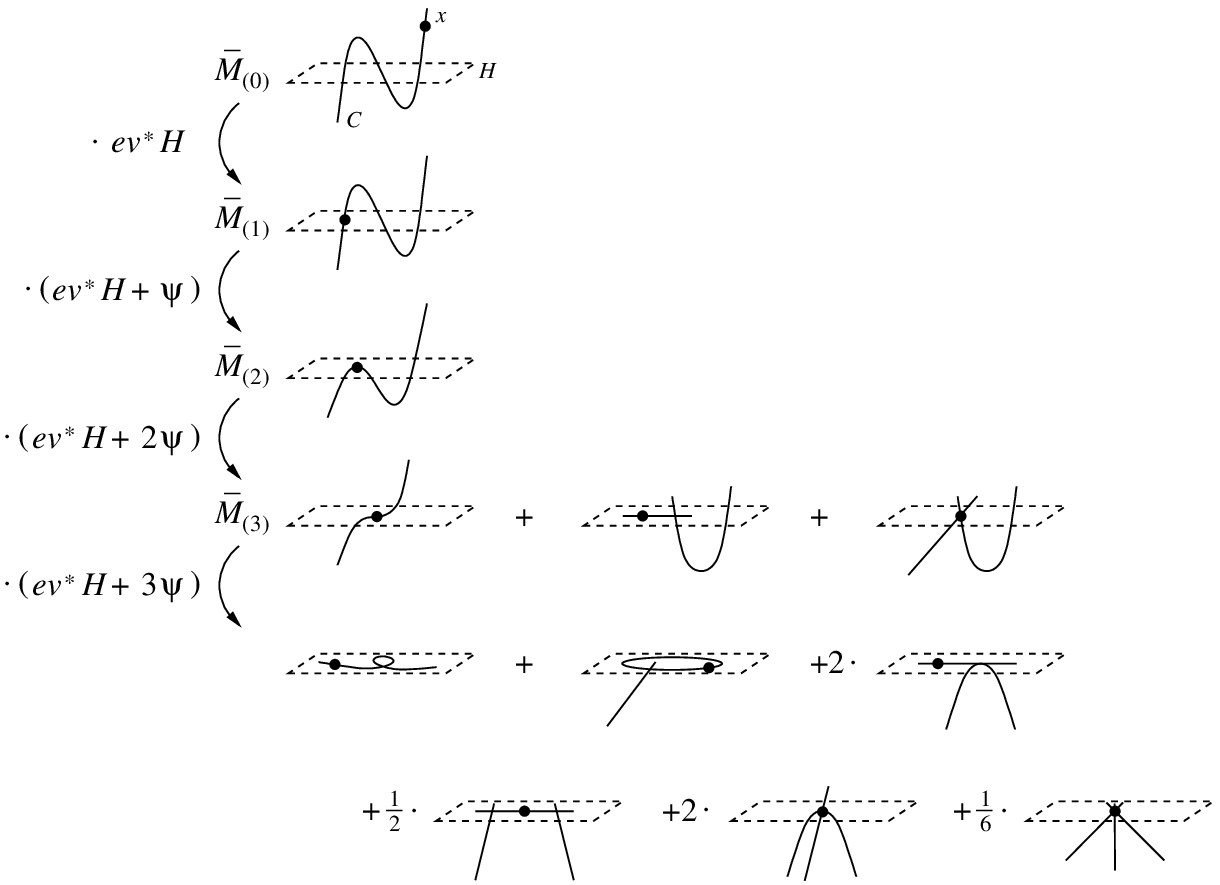} \end {center}
\vspace {0mm plus2mm}

(Of course, in the pictures where we have drawn the marked point on a node of
the curve, the corresponding stable maps have a contracted component, i.e.\ we
have $ \beta^{(0)} = 0 $.)

So we see that $ \bar M_1 (H,3) $ is equal to $ \prod_{i=0}^3 (ev^* H + i\psi)
· \bar M_1 (\PP^3,3) $ plus a bunch of correction terms coming from reducible
curves as shown in the picture. This is an equation of 9-dimensional cycles in
$ \bar M_1 (\PP^3,3) $. To make this into equations for the Gromov-Witten
invariants of $H$, we have to intersect it with some cohomology class $ \gamma
$ of codimension 9 that is a polynomial in $ ev^* H $ and $ \psi $. Note that
in the correction terms this will impose 9 conditions on the component $
C^{(0)} $ contained in $H$. However, in all the terms where the degree of $
C^{(0)} $ is at most 2, the moduli space for this component has dimension
smaller than 9. Hence all these terms vanish, and it follows that the 1-point
Gromov-Witten invariants of $H$ (of degree 3 in this example) are expressible
in terms of those of $ \PP^3 $ as
  \[ I^H_3 (\gamma) = I^{\PP^3}_3 \left(
       \gamma · \prod_{i=0}^3 (H + i\psi)
     \right). \]
The same argument works for higher degree of the curves.

Now let us come back to the case of general $X$ and $Y$. Can we still hope that
the correction terms vanish when we compute the Gromov-Witten invariants?
Recall that the reason for the vanishing above was that the dimension of the
moduli space of curves in $Y$ quickly gets bigger when the degree of the curves
goes up (in the example, the 9 conditions that were needed for Gromov-Witten
invariants for cubics in $Y$ were ``too many'' for lines and conics in $Y$).
Hence, as the (virtual) dimension of the moduli space of stable maps to $Y$ is
$ \vdim \bar M_1 (Y,\beta) = - K_Y·\beta + \dim Y - 2 $, we see that we need
that $ -K_Y $ is sufficiently positive.

If $ -K_Y $ is negative, basically all correction terms that could appear in
the computation of the Gromov-Witten invariants will do so. The main nuisance
about this is that the correction terms contain the full $n$-point
Gromov-Witten invariants of $Y$ (namely, $ n=1+r $ in each of the correction
terms), and not just the 1-point invariants that we originally wanted to
compute. There would be two ways to proceed:
\begin {itemize}
\item Use the version of theorem \ref {old-thm} for $n$-point invariants as
  proven in \cite {Ga}.
\item Use the WDVV equations to compute the $n$-point invariants of $Y$ in
  terms of 1-point invariants whenever they occur.
\end {itemize}
Both methods can be used without problems to write down an algorithm to compute
the Gromov-Witten invariants of $Y$ in terms of those of $X$. However, we do
not know at the moment how to express the result in a nice closed form.

Most interesting are the cases where $ -K_Y $ is nef, but yet not ``positive
enough'' to ensure the vanishing of all correction terms. We will show that,
whenever $ -K_Y $ is nef, the only $n$-point invariants of $Y$ that might occur
in the algorithm are those with fundamental or divisor classes at all but the
first marked point. These invariants can of course be reduced immediately to
1-point invariants using the fundamental class and divisor axioms for
Gromov-Witten invariants. Thus we arrive at recursion formulas that involve
only 1-point invariants. Solving them directly, we obtain a nice expression for
the invariants of $Y$: the ``mirror formula''.

The necessary computations to achieve this are done in section \ref
{sec-computations}. In section \ref {sec-examples} we apply the results to two
examples. First of all we rederive the expression for the genus zero
Gromov-Witten invariants of the quintic threefold. Secondly, we prove a similar
expression for the (virtual) numbers of plane rational curves of degree $d$
having contact of order $ 3d $ to a smooth cubic. These numbers play a role in
local mirror symmetry (see \cite {CKYZ} and \cite {T}). They are a by-product
of our work, as they are just simple examples of relative Gromov-Witten
invariants. The two main computational lemmas (that have nothing to do with
algebraic geometry, but rather are formal statements about certain power
series occurring in the calculation) are proved in the appendix.

The author would like to thank T. Graber, J. Harris, and R. Vakil for numerous
discussions. The work has been done at the Harvard University, to which the
author is grateful for hospitality.

  % Main file: latex main

\section {The mirror transformation} \label {sec-computations}

As in the introduction let $X$ be a smooth complex projective variety, and let
$Y$ be a smooth very ample hypersurface such that $ -K_Y $ is nef. By abuse of
notation, we denote by $ H^* (X) $ and $ H_* (X) $ the groups of algebraic
(co-)homology classes modulo torsion. For a class $ \beta \in H_2 (X) $ we
write $ \beta \ge 0 $ if $ \beta $ is effective, and $ \beta > 0 $ if $ \beta
\ge 0 $ and $ \beta \neq 0 $. To keep the notation as simple as possible, we
will assume in the following computations that the class of $Y$ generates $ H^2
(X) $ over $ \QQ $ (see remark \ref {general-h2} for the changes needed in the
general case).

For any $ \beta > 0 $ we denote by $ \bar M_n (X,\beta) $ the space of
$n$-pointed rational stable maps of class $ \beta $ to $X$. Let $ ev_i:
\bar M_n (X,\beta) \to X $ be the evaluation maps, and let $ \psi_i $ be the
cotangent line classes. For cohomology classes $ \gamma_i \in H^*(X) $ the
corresponding Gromov-Witten invariant is defined to be
  \[ I^X_\beta (\gamma_1 \psi^{k_1} \otimes \cdots \otimes \gamma_n \psi^{k_n})
       := ev_1^* \gamma_1 · \psi_1^{k_1} \cdots
          ev_n^* \gamma_n · \psi_n^{k_n} · [\bar M_n (X,\beta)]^{virt}
       \in \QQ \]
if the dimension condition $ \sum_i ( \codim \gamma_i + k_i ) = \vdim \bar M_n
(X,\beta) $ is satisfied, and zero otherwise. It is usual and convenient to
encode all the 1-point invariants of class $ \beta $ in a single cohomology
class
\begin {align*}
  I^X_\beta :=& \, ev_* \left( \frac 1 {1-\psi} ·
                [\bar M_1 (X,\beta)]^{virt} \right) \\
             =& \sum_{i,j} I^X_\beta (T^i \psi^j) · T_i
                \qquad \qquad \qquad \in H^* (X),
\end {align*}
where $ ev = ev_1 $, $ \{ T^i \} $ is a basis of $ H^*(X) \otimes \QQ $, and $
\{ T_i \} $ is the dual basis. Note that the dimension condition ensures that
for each $i$ at most one $j$ contributes a non-zero term to the sum above, so
all 1-point invariants of $X$ of class $ \beta $ can be reconstructed from the
cohomology class $ I^X_\beta $.

We define the Gromov-Witten invariants $ I^Y_\beta $ of $Y$ in the same way,
replacing $ \bar M_n (X,\beta) $ by $ \bar M_n (Y,\beta) $, \emph {but keeping
the $ ev_i $ to denote the evaluation maps to $X$}. Note that $ \beta $ is
still a homology class in $X$; so strictly speaking $ \bar M_n (Y,\beta) $ is
the space of stable maps to $Y$ of all homology classes whose push-forward to
$X$ is $ \beta $.

For $ \beta=0 $, we set $ I^X_0 := 1 $ and $ I^Y_0 := Y $.

Now consider the moduli spaces $ \bar M_{(m)} (X,\beta) $ of 1-pointed relative
stable maps to $X$ with multiplicity $m$ to $Y$ at the marked point (\cite {Ga}
definition 1.1). In the same manner as above, these spaces together with their
virtual fundamental classes (\cite {Ga} definition 1.18) give rise to
invariants $ I_{\beta,(m)} (\gamma \psi^k) $ that can be assembled into a
cohomology class
  \[ I_{\beta,(m)} = ev_* \left( \frac 1 {1-\psi} ·
                     [\bar M_{(m)} (X,\beta)]^{virt} \right)
                     \quad \in H^*(X). \]

\begin {remark} \label {vanish}
  For future reference, let us note that (as expected from geometry)
  $ I_{\beta,(0)} = I^X_\beta $ and $ I_{\beta,(m)} = 0 $ for $ m > Y·\beta $
  (see \cite {Ga} remark 1.3).
\end {remark}

Finally, let $ D_{(m)} (X,\beta) $ be the correction terms defined in theorem
\ref {old-thm}, and set
\begin {gather} \label {def-j}
  J_{\beta,(m)} = ev_* \left( \frac 1 {1-\psi} ·
                  [D_{(m)} (X,\beta)]^{virt} \right)
                + m · ev_* [\bar M_{(m)} (X,\beta)]^{virt} \quad \in H^*(X).
\end {gather}
The surprising additional term will appear in the proof of the following lemma.
Geometrically, it corresponds to unstable maps that have two irreducible
components $ C^{(0)} $ and $ C^{(1)} $, where $ C^{(0)} $ is contracted to a
point in $Y$ and contains the marked point, and $ C^{(1)} $ is a curve with
multiplicity $m$ to $Y$ at this point (see the end of the proof of lemma \ref
{j-comp}).

The first thing to do is to rewrite theorem \ref {old-thm} in the new
simplified notation.
\begin {lemma} \label {old-thm-new}
  For all $ \beta>0 $ and $ m \ge 0 $ we have
    \[ (Y+m) · I_{\beta,(m)} = I_{\beta,(m+1)} + J_{\beta,(m)}
       \quad \in H^*(X). \]
\end {lemma}

\begin {proof}
  Intersect the equation of theorem \ref {old-thm} with $ \frac 1 {1-\psi} $
  and push it forward by the evaluation map to get
  \begin {align*}
    ev_* &\left( (ev^* Y + m\psi) ·
      \frac 1 {1-\psi} · [\bar M_{(m)} (X,\beta)]^{virt} \right) \\
    &= ev_* \left( \frac 1 {1-\psi} · [\bar M_{(m+1)} (X,\beta)]^{virt} \right)
      + ev_* \left( \frac 1 {1-\psi} · [D_{(m)} (X,\beta)]^{virt} \right).
  \end {align*}
  As $ \frac \psi {1-\psi} = \frac 1 {1-\psi} - 1 $, the left hand side of this
  equation can be rewritten as
    \[ (Y+m) ·
         ev_* \left( \frac 1 {1-\psi} · [\bar M_{(m)} (X,\beta)]^{virt} \right)
       - m · ev_* [\bar M_{(m)} (X,\beta)]^{virt}. \]
  Taking into account the definitions of $ I_{\beta,(m)} $ and $ J_{\beta,(m)}
  $, we arrive at the equation stated in the lemma.
\end {proof}

\begin {remark} \label {thm-rec}
  In particular,
    \[ \prod_{i=0}^{Y·\beta} (Y+i) · I^X_\beta =
       \sum_{m=0}^{Y·\beta} \; \prod_{i=m+1}^{Y·\beta} (Y+i) · J_{\beta,(m)}.
    \]
  This follows from a recursive application of lemma \ref {old-thm-new},
  with the start and the end of the recursion given by remark \ref {vanish}.
\end {remark}

The next thing to do is to evaluate the $ J_{\beta,(m)} $ explicitly.

\begin {remark} \label {first-summand}
  Let us first consider the first summand $ ev_* \left( \frac 1 {1-\psi} ·
  [D_{(m)} (X,\beta)]^{virt} \right) $ in the definition (\ref {def-j}) of
  $ J_{\beta,(m)} $. Using the definition of $ D_{(m)} (X,\beta) $ and its
  virtual fundamental class given in theorem \ref {old-thm}, we see that this
  first summand is a sum of individual terms, each of which has the form
  \begin {gather} \label {sum-term}
    I^Y_{\beta^{(0)}} (T^i \psi^j \otimes \gamma_1 \otimes \cdots
      \otimes \gamma_r) · \frac 1 {r!} \; \prod_{k=1}^r \left(
        m^{(k)} · I_{\beta^{(k)},(m^{(k)})} (\gamma_k^\vee)
      \right) · T_i,
  \end {gather}
  where $ \gamma^\vee $ denotes the dual of a class $ \gamma $ \emph {in $Y$}.
  These terms are summed over all $i$, $ j \ge 0 $, $ r \ge 0 $, $ \beta^{(k)}
  $ (with $ \beta^{(0)} \ge 0 $ and $ \beta^{(k)} > 0 $ if $ k>0 $), and $
  m^{(k)} > 0 $, subject to the conditions
  \begin {enumerate}
  \item $ \beta^{(0)} + \cdots + \beta^{(r)} = \beta $ (degree condition),
  \item $ Y·\beta^{(0)} + m^{(1)} + \cdots + m^{(r)} = m $ (multiplicity
    condition),
  \item if $ \beta^{(0)} = 0 $ then $ r \ge 2 $ (stability condition).
  \end {enumerate}
  Moreover, the $ \gamma_k $ have to run over a basis of $ H^*(Y) \otimes \QQ $
  (actually it is sufficient to let them run over a basis of the part of $
  H^*(Y) \otimes \QQ $ induced by $X$, see \cite {Ga} remark 5.4).
\end {remark}

The main simplification of this huge sum is due to the following lemma, which
follows from a simple dimension count. It is the only point in our computations
where we need that $ - K_Y $ is nef.

\begin {lemma} \label {only-fund-div}
  The above expression (\ref {sum-term}) can only be non-zero if all $
  \gamma_k $ are fundamental or divisor classes. Moreover, for all $k$ we
  must have
    \[ \begin {array}{ll}
         m^{(k)} = Y·\beta^{(k)} - K_Y·\beta^{(k)} - 1 &
           \quad \mbox {if $ \gamma_k $ is the fundamental class}, \\
         m^{(k)} = Y·\beta^{(k)} - K_Y·\beta^{(k)} &
           \quad \mbox {if $ \gamma_k $ is a divisor class}.
       \end {array}
    \]
\end {lemma}

\begin {proof}
  As the invariants $ I_{\beta^{(k)},(m^{(k)})} (\gamma_k^\vee) $ must have
  dimension zero for all $k$, it follows that
  \begin {align*}
    \codim \gamma_k &= \dim Y - \codim \gamma_k^\vee \\
      &= \dim Y - \dim \bar M_{(m^{(k)})} (X,\beta^{(k)}) \\
      &= \dim Y - (- K_X·\beta^{(k)} + \dim X - 2 - m^{(k)}) \\
      &= - Y·\beta^{(k)} + K_Y·\beta^{(k)} + 1 + m^{(k)}
          \qquad \mbox {(by adjunction)}.
  \end {align*}
  This shows the equation for the $ m^{(k)} $. Moreover, as $ -K_Y $ is nef and
  we must have $ m^{(k)} \le Y·\beta^{(k)} $ for the relative invariant to be
  non-zero (see remark \ref {vanish}), it follows that $ \codim \gamma_k \le 1
  $, as desired.
\end {proof}

\begin {remark} \label {special}
  Obviously, in the same way one can show that:
  \begin {itemize}
  \item If $ - K_Y·\beta \ge 1 $ for all $ \beta > 0 $ then all the $ \gamma_k
    $ have to be fundamental classes. (In the following computations this
    would mean that all $ r_\beta = 0 $, which greatly simplifies the
    calculation.) This is e.g.\ the case if $Y$ is a hypersurface in $ X=\PP^n
    $ of degree at most $n$.
  \item If $ - K_Y·\beta \ge 2 $ for all $ \beta > 0 $ then no $ \gamma_k $
    can exist, i.e.\ we must always have $ r=0 $. Hence in this case we
    conclude that there are no correction terms in the computation of the
    Gromov-Witten invariants. The only term on the right hand side of remark
    \ref {thm-rec} is $ I^Y_\beta $ (for $ r=0 $ and $ m=Y·\beta $), so it
    follows that the ``naïve'' formula
      \[ I^Y_\beta = \prod_{i=0}^{Y·\beta} (Y+i) · I^X_\beta \]
    is true (as in the case considered in the introduction where $ Y \subset X
    $ is a plane in $ \PP^3 $). This is e.g.\ the case if $Y$ is a hypersurface
    in $ X=\PP^n $ of degree at most $ n-1 $.
  \end {itemize}
\end {remark}

\begin {remark} \label {multi-index}
  As we have assumed that the class of $Y$ generates $ H^2 (X) $ over $ \QQ $,
  lemma \ref {only-fund-div} states that the only factors that can occur in the
  $k$-product in (\ref {sum-term}) are the numbers
  \begin {align*}
    s_\beta &:= (Y·\beta - K_Y·\beta - 1) ·
                I_{\beta,(Y·\beta - K_Y·\beta - 1)} (1^\vee) \\
    \mbox {and} \quad
    r_\beta &:= (Y·\beta - K_Y·\beta) ·
                I_{\beta,(Y·\beta - K_Y·\beta)} (Y^\vee)
  \end {align*}
  for some $ \beta > 0 $. Thus we can then rewrite (\ref {sum-term}) using
  multi-index notation as follows. For a multi-index $ \mu = (\mu_\beta) $ of
  non-negative integers indexed by the positive homology classes $ \beta $ of
  $ H^2(X) $, we apply the usual notations
    \[ \begin {array}{r@{\;}c@{\;}l@{\qquad \qquad}r@{\;}c@{\;}l}
    \sum \mu &:=& \sum_\beta \mu_\beta, &
    s^\mu &:=& \prod_\beta s_\beta^{\mu_\beta}, \\
    \mu! &:=& \prod_\beta \mu_\beta!, &
    |\mu| &:=& \sum_\beta \mu_\beta · \beta.
       \end {array} \]
  Then we can rewrite (\ref {sum-term}) as
  \begin {gather} \label {alt-term}
    I^Y_{\beta^{(0)}} (T^i \psi^j \otimes 1^{\otimes \sum \mu}
      \otimes Y^{\otimes \sum \nu}) · \frac 1 {r!} ·
      s^\mu r^\nu · T_i,
  \end {gather}
  where $ \mu $ and $ \nu $ are the multi-indices such that the factors $
  s_\beta $ and $ r_\beta $ appear in (\ref {sum-term}) $ \mu_\beta $ and $
  \nu_\beta $ times, respectively. In particular, $ r = \sum \mu + \sum \nu $
  is the number of nodes of the curves under consideration.
\end {remark}

We are now ready to evaluate the $ J_{\beta,(m)} $ explicitly in terms of the
1-point Gromov-Witten invariants $ I^Y_\beta $ of $Y$ and the relative 1-point
invariants $ s_\beta $ and $ r_\beta $.

\begin {lemma} \label {j-comp}
  With the notation of remark \ref {multi-index},
    \[ J_{\beta,(m)} = \sum_{\mu,\nu}
                       \left( Y + Y·\beta^{(0)} \right)^{\sum \nu}
                     · \frac {s^\mu}{\mu!} \frac {r^\nu}{\nu!}
                     · I^Y_{\beta^{(0)}} \]
  for all $ \beta > 0 $ and $ m \ge 0 $, where the sum is taken over all
  multi-indices $ \mu $ and $ \nu $ such that $ \beta^{(0)} := \beta - |\mu| -
  |\nu| \ge 0 $ (degree condition) and $ m = Y·\beta - K_Y·(|\mu|+|\nu|) -
  \sum \mu $ (multiplicity condition).
\end {lemma}

\begin {proof}
  Inserting expression (\ref {alt-term}) for (\ref {sum-term}) in remark \ref
  {first-summand}, we see that the first summand in the definition (\ref
  {def-j}) of $ J_{\beta,(m)} $ is
  \begin {align*}
    & ev_* \left( \frac 1 {1-\psi} ·
      [D_{(m)} (X,\beta)]^{virt} \right) \\
    & \qquad \qquad \qquad
    = \sum_{i,j} \sum_{\mu,\nu}
      I^Y_{\beta^{(0)}} (T^i \psi^j \otimes
      1^{\otimes \sum \mu} \otimes Y^{\otimes \sum \nu})
    · \frac {s^\mu}{\mu!} \frac {r^\nu}{\nu!}
    · T_i,
  \end {align*}
  where the sum is taken over all $ i,j,\mu,\nu $ such that
  \begin {enumerate}
  \item $ \beta^{(0)} := \beta - |\mu| - |\nu| \ge 0 $ (degree condition),
  \item $ Y·\beta - K_Y·(|\mu|+|\nu|) - \sum \mu = m $ (multiplicity
    condition --- here we inserted the expression of lemma \ref {only-fund-div}
    for the $ m^{(i)} $),
  \item if $ \beta^{(0)} = 0 $ then $ \sum \mu + \sum \nu \ge 2 $ (stability
    condition).
  \end {enumerate}
  Now we compute the Gromov-Witten invariant $ I^Y_{\beta^{(0)}} (\cdots) $ in
  terms of 1-point invariants of $Y$. We claim that for $ \beta^{(0)} > 0 $
  \begin {gather} \label {div-axiom}
    \sum_{i,j} I^Y_{\beta^{(0)}} (T^i \psi^j \otimes
      1^{\otimes \sum \mu} \otimes Y^{\otimes \sum \nu}) · T_i
    = ( Y+Y·\beta^{(0)} )^{\sum \nu} · I^Y_{\beta^{(0)}}.
  \end {gather}
  In fact, this follows from the fundamental class axiom
  \begin {align*}
    \sum_{i,j} I^Y_{\beta^{(0)}} (T^i \psi^j \otimes 1 \otimes \cdots) · T_i
    &= \sum_{i,j \neq 0} I^Y_{\beta^{(0)}} (T^i \psi^{j-1} \otimes \cdots)
       · T_i \\
    &= \sum_{i,j} I^Y_{\beta^{(0)}} (T^i \psi^j \otimes \cdots) · T_i
  \end {align*}
  and the divisor axiom
  \begin {align*}
    \sum_{i,j} I^Y_{\beta^{(0)}} (T^i \psi^j \otimes Y \otimes \cdots) · T_i
    &= \sum_{i,j} (Y·\beta^{(0)}) ·
         I^Y_{\beta^{(0)}} (T^i \psi^j \otimes \cdots) · T_i \\
    &\qquad + \sum_{i,j \neq 0} I^Y_{\beta^{(0)}} (T^i · Y \; \psi^{j-1}
         \otimes \cdots) · T_i \\
    &= \sum_{i,j} (Y·\beta^{(0)}) ·
         I^Y_{\beta^{(0)}} (T^i \psi^j \otimes \cdots) · T_i \\
    &\qquad + \sum_{i,j \neq 0} I^Y_{\beta^{(0)}} (T^i \psi^{j-1}
         \otimes \cdots) · (T_i·Y) \\
    &= (Y·\beta^{(0)} + Y) ·
       \sum_{i,j} I^Y_{\beta^{(0)}} (T^i \psi^j \otimes \cdots)·T_i
  \end {align*}
  (see e.g.\ \cite {Ge} proposition 12), where ``$\cdots$'' denotes any tensor
  product of cohomology classes (i.e.\ not including cotangent line classes).
  In fact, the same formula (\ref {div-axiom}) is also true for $ \beta^{(0)}
  = 0 $, as in this case
  \begin {align*}
    \sum_{i,j} I^Y_0 (T^i \psi^j \otimes
       1^{\otimes \sum \mu} \otimes Y^{\otimes \sum \nu}) · T_i
    &= (Y^{\sum \nu}) · Y \\
    &= Y^{\sum \nu} · I^Y_0
  \end {align*}
  by the ``mapping to point axiom''. Hence the first summand in the definition
  (\ref {def-j}) of $ J_{\beta,(m)} $ is
  \begin {gather} \label {sum-first}
    ev_* \left( \frac 1 {1-\psi} ·
     [D_{(m)} (X,\beta)]^{virt} \right)
    = \sum_{\mu,\nu}
      \left( Y + Y·\beta^{(0)} \right)^{\sum \nu}
      · \frac {s^\mu}{\mu!} \frac {r^\nu}{\nu!}
      · I^Y_{\beta^{(0)}}
  \end {gather}
  with the sum taken over all $ \mu,\nu $ satisfying the degree, multiplicity
  and stability conditions. The second summand is
  \begin {align*}
    m · ev_* [\bar M_{(m)} (X,\beta)]^{virt}
      &= m · \sum_i I_{\beta,(m)} (T^i) · T_i \\
      &= s_\beta · Y · \delta_{m,Y·\beta-K_Y·\beta-1} \\
      &\qquad + r_\beta · Y^2 · \delta_{m,Y·\beta-K_Y·\beta}
  \end {align*}
  by lemma \ref {only-fund-div}. As we have defined $ I^Y_0 = Y $, this adds
  exactly the terms with $ \beta^{(0)} = 0 $ and $ \sum \mu + \sum \nu = 1 $ to
  the sum in (\ref {sum-first}) that were excluded because of the stability
  condition. It follows that
    \[ J_{\beta,(m)} = \sum_{\mu,\nu}
                       \left( Y + Y·\beta^{(0)} \right)^{\sum \nu}
                     · \frac {s^\mu}{\mu!} \frac {r^\nu}{\nu!}
                     · I^Y_{\beta^{(0)}}, \]
  with the sum taken over all $ \mu,\nu $ satisfying the degree and
  multiplicity conditions.
\end {proof}

\begin {remark} \label {alt-mult}
  The multiplicity condition in lemma \ref {j-comp} can be replaced by
    \[ m = Y·\beta - \epsilon \sum \mu, \]
  where $ \epsilon \in \{0,1\} $ depends only on $Y$. To see this, recall that
  the multiplicity condition was obtained from the original one
  \begin {gather} \label {oldmult}
    m = Y·\beta^{(0)} + \sum m^{(k)}
  \end {gather}
  by inserting the expressions $ m^{(k)} = Y·\beta^{(k)} - K_Y·\beta^{(k)} $
  (for every $ r_{\beta^{(k)}} $) or $ m^{(k)} = Y·\beta^{(k)} -
  K_Y·\beta^{(k)} -1 $ (for every $ s_{\beta^{(k)}} $), respectively. But by
  remark \ref {vanish} we have $ r_{\beta^{(k)}} = 0 $ if $ m^{(k)} =
  Y·\beta^{(k)} - K_Y·\beta^{(k)} > Y·\beta^{(k)} $. So (as $ K_Y $ is nef) $
  r_{\beta^{(k)}} $ can only be non-zero if $ m^{(k)} = Y·\beta^{(k)} $. Hence
  we can insert this simplified expression for $ m^{(k)} $ in (\ref {oldmult}).

  In the same way, $ s_{\beta^{(k)}} $ can only be non-zero if $ m^{(k)} =
  Y·\beta^{(k)} -1 $ (in the case $ K_Y = 0 $) or $ m^{(k)} = Y·\beta^{(k)} $
  (in the case $ K_Y > 0 $). In other words, $ m^{(k)} = Y·\beta^{(k)} -
  \epsilon $ with $ \epsilon \in \{0,1\} $ depending only on $Y$.
  
  If we now take the original multiplicity condition (\ref {oldmult}) and
  insert the new simplified expressions $ m^{(k)} = Y·\beta^{(k)} $ (for every
  $ r_{\beta^{(k)}} $) and $ m^{(k)} = Y·\beta^{(k)} - \epsilon $ (for every $
  s_{\beta^{(k)}} $), respectively, we arrive at the desired multiplicity
  condition $ m = Y·\beta - \epsilon \sum \mu $.
\end {remark}

\begin {remark} \label {insert}
  Now we can insert the expression of lemma \ref {j-comp} (with the
  multiplicity condition from remark \ref {alt-mult}) into the formula of
  remark \ref {thm-rec}. Thus we obtain
  \begin {align*}
    \prod_{i=0}^{Y·\beta} (Y+i) · I^X_\beta
      &= \sum_{\mu,\nu}
         \prod_{i=Y·\beta - \epsilon \sum \mu + 1}^{Y·\beta} (Y+i)
       · \left( Y + Y·\beta^{(0)} \right)^{\sum \nu}
       · \frac {s^\mu}{\mu!} \frac {r^\nu}{\nu!}
       · I^Y_{\beta^{(0)}} \\
      &= \sum_{\mu,\nu}
         \prod_{i=0}^{\epsilon \sum \mu - 1} ( Y + Y·\beta - i )
       · \left( Y + Y·\beta^{(0)} \right)^{\sum \nu}
       · \frac {s^\mu}{\mu!} \frac {r^\nu}{\nu!}
       · I^Y_{\beta^{(0)}},
  \end {align*}
  where the sum is now taken over all $ \mu,\nu $ satisfying the degree
  condition $ \beta^{(0)} := \beta - |\mu|-|\nu| \ge 0 $. Note that this
  equation is trivially true in the case $ \beta = 0 $ as well (both sides are
  equal to $Y$ in this case).

  To get rid of the degree condition, we multiply these equations with $
  q^{Y·\beta} $ (where $q$ is a formal variable) and add them up; so we get
  \begin {align} \label {thm-final}
    \sum_\beta & \prod_{i=0}^{Y·\beta} (Y+i) · I^X_\beta · q^{Y·\beta}
         \notag \\
      &= \sum_{\beta^{(0)}} \sum_{\mu,\nu}
         \prod_{i=0}^{\epsilon \sum \mu - 1} ( Y + Y·\beta - i )
       · \left( Y + Y·\beta^{(0)} \right)^{\sum \nu}
       · \frac {s^\mu}{\mu!} \frac {r^\nu}{\nu!}
       · I^Y_{\beta^{(0)}}
       · q^{Y·\beta},
  \end {align}
  where the sum now runs over all multi-indices $ \mu,\nu $ (and $ \beta =
  \beta^{(0)} +|\mu|+|\nu| $).
  
  Although this equation looks quite complicated, note that all geometric ideas
  in its derivation are still visible: the left hand side is the ``naïve''
  expression for the Gromov-Witten invariants of $Y$ that we already
  encountered in the introduction and remark \ref {special}. The product $
  \prod_{i=0}^{Y·\beta} (Y+i) $ here corresponds to the process of raising the
  multiplicity of the curves from $0$ to $ Y·\beta+1 $. The right hand side of
  the equation describes the correction terms. They correspond to reducible
  curves with one component in the hypersurface ($ I^Y_{\beta^{(0)}} $) and
  various others in the ambient space with specified multiplicities to the
  hypersurface ($ s^\mu r^\nu $). The factor $ (Y+Y·\beta^{(0)})^{\sum \nu} $
  comes from the $ (\sum \nu) $-fold application of the divisor axiom that we
  used to describe the component in the hypersurface by a 1-point invariant
  instead of by a $ (1+r) $-point invariant.
  
  All that remains to be done to arrive at the ``mirror formula'' is to
  simplify the right hand side of equation (\ref {thm-final}). To do so, define
  $ P(t) $ to be ``the right hand side with $ Y·\beta^{(0)} $ replaced by a
  formal variable $t$'':
\end {remark}

\begin {definition} \label {def-p}
  Let
    \[ P(t) := \sum_{\mu,\nu}
               \prod_{i=0}^{\epsilon \sum \mu - 1}
                 ( Y + Y·(|\mu|+|\nu|) + t - i )
             · \left( Y + t \right)^{\sum \nu}
             · \frac {s^\mu}{\mu!} \frac {r^\nu}{\nu!}
             · q^{Y·(|\mu|+|\nu|)}, \]
  so that (\ref {thm-final}) can be written as
  \begin {gather} \label {thm-p}
    \sum_\beta \prod_{i=0}^{Y·\beta} (Y+i) · I^X_\beta · q^{Y·\beta}
      = \sum_{\beta} P(Y·\beta)
      · I^Y_\beta · q^{Y·\beta}.
  \end {gather}
\end {definition}

\begin {lemma} \label {diff-eqn}
  The power series $ P(t) $ of definition \ref {def-p} satisfies the
  differential equation $ \frac {d^2}{dt^2} \ln P = 0 $. In particular, if $
  P(t) = P_0 + P_1·t + \cdots $ is the Taylor expansion of $P$ then $ P(t) =
  P_0 \exp ( \frac {P_1}{P_0} t) $.
\end {lemma}

\begin {proof}
  This can be checked directly from the definition of $ P(t) $. The statement
  does not depend on the special values of $ r_\beta $ and $ s_\beta $; it is
  equally true if the $ r_\beta $ and $ s_\beta $ are considered to be formal
  variables. We give a proof of the statement in appendix \ref {sec-lemma}
  (apply lemma \ref {tech-lemma} with the collection of variables $ x_i $ being
  the union of the $ r_\beta $ and $ s_\beta $, $ z=0 $, and $t$ replaced by $
  t+Y $).
\end {proof}

\begin {corollary}[Mirror formula] \label {mirror-cor}
  If we formally set $ \tilde q = q·\exp \frac {P_1}{P_0} $ with $ P_0 $ and $
  P_1 $ as in lemma \ref {diff-eqn}, then
    \[ \sum_\beta \prod_{i=0}^{Y·\beta} (Y+i) · I^X_\beta · q^{Y·\beta}
       = P_0 · \sum_{\beta} I^Y_\beta · \tilde q^{Y·\beta}, \]
  i.e.\ the generating function $ \sum_{\beta} I^Y_\beta · q^{Y·\beta} $ of the
  1-point Gromov-Witten invariants of $Y$ can be obtained from the ``naïve''
  expression $ \sum_\beta \prod_{i=0}^{Y·\beta} (Y+i) · I^X_\beta · q^{Y·\beta}
  $ by a formal change of variables ($ q \to \tilde q $) and a scaling factor
  ($ \cdot P_0 $).
\end {corollary}

\begin {proof}
  Immediately from (\ref {thm-p}) and lemma \ref {diff-eqn}.
\end {proof}

\begin {remark} \label {general-h2}
  In the above computations we assumed that the class of $Y$ generates $ H^2
  (X) $ over $ \QQ $. In fact, this is not essential. All that happens for
  higher dimension of $ H^2(X) $ is that the notation becomes more complicated
  at some steps of the calculation. Most importantly, in remark \ref
  {multi-index} there are now more factors that can occur in the $k$-product of
  (\ref {sum-term}). Namely, instead of the $ r_\beta $ we now have
      \[ r_{i,\beta} = (Y·\beta - K_Y·\beta) ·
                        I_{\beta,(Y·\beta - K_Y·\beta)} (\gamma_i^\vee), \]
  for $ i=1,\dots,\dim H^2(X) \otimes \QQ $, where the $ \gamma_i $ form a
  basis of $ H^2(X) \otimes \QQ $, chosen such that $ \gamma_1 = Y $.
  Correspondingly, lemma \ref {j-comp} becomes
    \[ J_{\beta,(m)} = \sum_{\mu,\nu_i} \prod_i
                       \left( \gamma_i + \gamma_i·\beta^{(0)} \right)^{
                         \sum \nu_i}
                     · \frac {s^\mu}{\mu!}
                     · \prod_i \frac {{r_i}^{\nu_i}}{\nu_i!}
                     · I^Y_{\beta^{(0)}} \]
  where the $ \nu_i $ are multi-indices. In the alternative multiplicity
  condition of remark \ref {alt-mult}, the number $ \epsilon $ will now depend
  on $ \beta $ (it is 1 if $ K_Y·\beta = 0 $ and 0 if $ K_Y·\beta >0 $). Hence
  the multiplicity condition is now $ m = Y·\beta - \epsilon \mu $, where $
  \epsilon $ is a multi-index with entries 0 and 1. Finally, we need a formal
  variable $ q_i $ for each $ \gamma_i $ to replace the expression $
  q^{Y·\beta} $ by $ q^\beta := \prod_i q_i^{\gamma_i·\beta} $.
  Definition \ref {def-p} then becomes
  \begin {align*}
    P(\{t_i\}) :=& \sum_{\mu,\nu_i}
                    \prod_{j=0}^{\epsilon \mu - 1}
                    ( Y + Y·(|\mu|+\sum_i |\nu_i|) + t_1 - j )
                · \prod_i \left( \gamma_i + t_i \right)^{\sum \nu_i} \\
                &\qquad · \frac {s^\mu}{\mu!}
                · \prod_i \frac {{r_i}^{\nu_i}}{\nu_i!}
                · q^{|\mu|+\sum_i |\nu_i|},
  \end {align*}
  with which we obtain the equation (compare to (\ref {thm-p}))
  \begin {gather}
    \sum_\beta \prod_{i=0}^{Y·\beta} (Y+i) · I^X_\beta · q^\beta
      = \sum_{\beta} P(\{\gamma_i·\beta\})
      · I^Y_\beta · q^\beta.
  \end {gather}
  The same proof as for lemma \ref {diff-eqn} works to show that $ \partial_{
  t_i} \partial_{t_j} \ln P = 0 $ for all $ i,j $, so it follows that $ P(t) =
  P_0 \exp (\frac {\sum P_i t_i}{P_0}) $, where $ P(\{t_i\}) = P_0 + \sum_i
  P_i·t_i + \cdots $ is the linear expansion of $P$. Hence the mirror formula
  of corollary \ref {mirror-cor} holds in the same way
    \[ \sum_\beta \prod_{i=0}^{Y·\beta} (Y+i) · I^X_\beta · q^\beta
       = P_0 · \sum_{\beta} I^Y_\beta · \tilde q^\beta, \]
  where $ \tilde q_i = q_i · \exp \frac {P_i}{P_0} $.
\end {remark}

  % Main file: latex main

\section {Examples} \label {sec-examples}

\begin {example}[Application to the quintic threefold] \label {quintic}
  Let $ X=\PP^4 $, and let $ Y \subset X $ be a smooth quintic hypersurface, so
  that $ Y = 5H \in H^*(X) $, where $H$ is the class of a hyperplane. We are
  interested in the genus zero Gromov-Witten invariants of $Y$, i.e.\ in the
  numbers $ n_d = \frac 1 d I^Y_d (H) $ (note that $H$ has $d$ points of
  intersection with a degree-$d$ curve). As this is the $ H^3 $-coefficient of
  $ I^Y_d $ (up to a scaling factor), we consider the equation (\ref {thm-p})
  modulo $ H^4 $. (This discards the invariants $ I^Y_d (\psi) $.)

  Since the only Gromov-Witten invariants of $Y$ are $ I^Y_d (H) $ (and $ I^Y_d
  (\psi) $), the polynomials $ I^Y_d $ have no $ H^0 $, $ H^1 $, and $ H^2 $
  terms for $ d>0 $. Hence as it is well-known that
    \[ I^X_d = \prod_{i=1}^d \frac 1 {(H+i)^5}, \]
  (see e.g.\ \cite {P} section 1.4) it follows from (\ref {thm-p}) that
    \[ \sum_{d \ge 0} 5H
       · \frac {\prod_{i=1}^{5d} (5H+i)}{\prod_{i=1}^d (H+i)^5} \; q^{5d}
         = 5H \; P_0 \pmod {H^3}. \]
  This is sufficient to reconstruct $P$: if we expand
  \begin {gather} \label {def-f}
    \sum_{d \ge 0}
        \frac {\prod_{i=1}^{5d} (5H+i)}{\prod_{i=1}^d (H+i)^5} \; q^{5d}
      =: F_0 + F_1 H + F_2 H^2 + \cdots
  \end {gather}
  then $ P|_{t=H=0} = F_0 $ and $ \partial_H P|_{t=H=0} = F_1 $. So as $P$ is a
  function of $ t+5H $ and satisfies $ \partial_t^2 \ln P = 0 $,
  it follows that $ \partial_t P|_{t=H=0} = \frac 1 5 \, F_1 $, and
  hence
    \[ P = F_0 · \exp \left( \left( \frac t 5 + H \right) · \frac {F_1}{F_0}
         \right). \]
  In particular,
  \begin {align*}
    P_0 &= F_0 · \exp \left( H \, \frac {F_1}{F_0} \right) \\
        &= F_0 + H \, F_1 + \frac {H^2} 2 \,
           \frac {F_1^2}{F_0} + \cdots.
  \end {align*}
  So by comparing the $ H^3 $-coefficient of (\ref {thm-p}) we get
    \[ F_2 = \frac 1 2 \frac {F_1^2}{F_0}
           + \frac 1 5 \sum_{d>0} d n_d q^{5d} F_0 \exp
           \left( d \, \frac {F_1}{F_0} \right). \]
  Together with (\ref {def-f}), this equation determines the $ n_d $
  recursively and gives the well-known numbers $ n_1 = 2875 $, $ n_2 = 609250
  + \frac {2875}{8} $,\dots .
\end {example}

\begin {example}[Application to plane elliptic curves] \label {elliptic}
  We want to compute the (virtual) numbers of rational plane curves of degree
  $d$ having multiplicity $ 3d $ to a smooth elliptic plane cubic, i.e.\ the
  relative Gromov-Witten invariants $ I_{d,(3d)} (1) = \frac 3 d \, r_d $ in
  the case where $ X=\PP^2 $ and $Y$ is a smooth elliptic cubic. According to
  \cite {T} remark 1.11 these numbers are related to the local mirror symmetry
  of \cite {CKYZ}.

  The computation of the numbers $ r_d $ is very similar (yet not identical) to
  that of the Gromov-Witten invariants of $Y$ in section \ref
  {sec-computations}. This time we apply lemma \ref {old-thm-new} recursively
  only up to multiplicity $ 3d $ instead of $ 3d+1 $, so we get
    \[ \prod_{i=0}^{3d-1} (3H + i) \; I^X_d =
         I_{d,(3d)} + 
         \sum_{m=0}^{3d-1} \prod_{i=m+1}^{3d-1}(3H + i) \; J_{d,(m)}. \]
  Note that $ I^Y_d = 0 $ for $ d>0 $, as there are no rational curves in $Y$.
  So if we insert the expression for $ J_{d,(m)} $ of lemma \ref {j-comp}, we
  get in the same way as in remark \ref {insert}
  \begin {align} \label {eqn-rel}
    \sum_{d>0} \prod_{i=0}^{3d-1} (3H + i) \; I^X_d \, q^{3d}
    &= \sum_{d>0} \frac {3H^2}{d} \, r_d q^{3d} \notag \\
    &\qquad +
      \sum_{\mu,\nu} \prod_{i=3d-\sum \mu+1}^{3d-1} (3H + i) \;
      (3H)^{\sum \nu}
      · \frac {s^\mu}{\mu!} \frac {r^\nu}{\nu!} · 3H \, q^{3d}
  \end {align}
  where we already inserted the expression $ m=3d-\sum \mu $ for Calabi-Yau
  hypersurfaces (see remark \ref {alt-mult}). Here, in the second line we
  set $ d = |\mu|+|\nu| $, and we obviously only sum over those $ \mu $ with $
  \sum \mu \ge 1 $.

  Similar to definition \ref {def-p} let us set
    \[ Q(t) := \sum_\mu \prod_{i=1}^{\sum \mu-1} (3H·|\mu|+t-i) \;
               \frac {s^\mu}{\mu!} \; q^{3H·|\mu|} \, t, \]
  where the sum is now taken over all $ \mu $ --- not only those with $ \sum
  \mu \ge 1 $. The $ \mu=0 $ term contributes a 1 (together with the factor
  $t$). The definition of $ Q(t) $ is so that $ Q(3H) - 1 $ yields exactly
  the $ \nu=0 $ terms in the second line of (\ref {eqn-rel}).

  Similarly to lemma \ref {diff-eqn} the power series $ Q(t) $ satisfies a
  differential equation: by lemma \ref {tech-lemma-2} $ \ln Q(t) $ is linear in
  $t$, i.e.\ $ Q(t) = \exp (c·t) $. To compute $c$, we expand as in example
  \ref {quintic} the left hand side of (\ref {eqn-rel})
    \[ \sum_{d>0} \; 3H ·
         \frac {\prod_{i=1}^{3d-1} (3H+i)}{\prod_{i=1}^d (H+i)^3} \; q^{3d}
         =: F_1 H + F_2 H^2 + \cdots \]
  (in \cite {T} $ F_1(q^3) $ is called $ I_2^{(0)}(z) $, and $ F_2(q^3) $ is
  called $ I_3^{(0)}(z) $). As the $t$-expansion of $ Q(t) $ is
    \[ Q(t) = 1 + c t + \frac 1 2 c^2 t^2 + \cdots, \]
  comparison of the $ H^1 $ terms in (\ref {eqn-rel}) gives $ F_1 = (\mbox
  {the $ H^1 $ term of $ Q(3H) $}) = 3c $; so $ Q(t) = \exp (\frac {F_1·t} 3)
  $.

  Now compare the $ H^2 $ term in (\ref {eqn-rel}). Note that we must have $
  \sum \nu \le 1 $ because of the factor $ (3H)^{\sum \nu + 1} $. The $ \nu=0 $
  term is exactly the second coefficient of $ Q(3H) $ as remarked above, i.e.\
  $ \frac 1 2 \, F_1^2 $. The terms with $ \sum \nu = 1 $ can be written as a
  sum over $d$, where $d$ is the index of the one non-zero entry of $ \nu $.
  The contribution for a given $d$ is exactly $ 9 r_d q^{3d} \frac {Q(3d)}{3d}
  = \frac 3 d r_d q^{3d} \exp (d \, F_1) $, with the $ \mu=0 $ term in $ Q(3d)
  $ coming from the right hand side of the first line of (\ref {eqn-rel}). Thus
  we get the equation
    \[ F_2 = \frac 1 2 \, F_1^2
             + \sum_{d>0} \frac 3 d \, r_d \, q^{3d} \exp (d F_1), \]
  which determines the numbers $ \frac 3 d r_d = I_{d,(3d)} (1) $. The first
  few numbers are given in the following table.
  \begin {center} \begin {tabular}{|l|c|c|c|c|c|c|c|c|} \hline
    $d$ & 1 & 2 & 3 & 4 & 5 & 6 & 7 & 8 \\ \hline
    $ I_{d,(3d)}(1) $ & $ 9 $ & $ \frac {135} 4 $ & $ 244 $ &
                    $ \frac {36999}{16} $ & $ \frac {635634}{25} $ &
                    $ 307095 $ & $ \frac {193919175}{49} $ &
                    $ \frac {3422490759}{64} $
                    \\ \hline
  \end {tabular} \end {center}
  This equation is equivalent to the conjecture of remark 1.11 in \cite {T}.
  Together with \cite {T} theorem 2.1 it proves that $ I_{d,(3d)} = (-1)^d 3d
  K_d $, where $ K_d $ is the top Chern class of the rank-$ (3d-1) $ bundle on
  $ \bar M_0 (\PP^2,d) $ with fiber $ H^1 (C,f^* K_{\PP^2}) $ at the point $
  (C,f) \in \bar M_0 (\PP^2,d) $. At the moment we do not know of a geometric
  proof of this statement.
\end {example}

  \begin {appendix}
  % Main file: latex main

\section {Proof of the main technical lemmas} \label {sec-lemma}

In this appendix we show that the power series $ P(t) $ and $ Q(t) $ of
definition \ref {def-p} and example \ref {elliptic} satisfy certain
differential equations.

\begin {lemma} \label {tech-lemma}
  Let $ x_i $ be a collection of variables (possibly infinite), and let $ a_i,
  b_i \in \NN $, $ c_i \in \CC $. Define
    \[ P(t,z) = \sum_k \frac {x^k}{k!} \; t^{ak} \;
                  \prod_{i=0}^{bk-1} (ck+z+t-i), \]
  where $k$ is a multi-index, and where we used the usual multi-index notations
  $ ak = \sum_i a_i k_i $, $ x^k = \prod_i x_i^{k_i} $, $ k! = \prod_i k_i! $.
  Assume that, for every $i$, the pair $ (a_i,b_i) $ is $ (0,0) $, $ (1,0) $,
  or $ (0,1) $. Then
    \[ \dt^2 \ln P = \dz^2 \ln P = \dt \dz \ln P = 0. \]
\end {lemma}

\begin {proof}
  \emph {Step 1.} We consider the $ c_i $ to be formal variables and show by
  induction on $n$ that for every $i$ and every $ n \ge 0 $
    \[ \begin {array}{lr@{\;}r@{\;}r@{\;}l}
    \mbox {if} & \dt^2 \ln P \at = & \dz^2 \ln P \at = &
        \dt \dz \ln P \at = & 0 \\
    \mbox {then} & \dci^n \dt^2 \ln P \at = & \dci^n \dz^2 \ln P \at = &
        \dci^n \dt \dz \ln P \at = & 0.
       \end {array} \]
  So assume that
    \[ \dci^j \dt^2 \ln P \at = \dci^j \dz^2 \ln P \at =
        \dci^j \dt \dz \ln P \at = 0 \]
  for $ j \le n $. Note that by definition of $P$ we have $ \dci P = x_i \dxi
  \dz P $. Let $ \di $ and $ \dj $ denote either $ \dt $ or $ \dz $. Then it
  follows that (everything in the following calculation is evaluated at $ c_i
  = 0 $):
  \begin {align*}
    \dci^{n+1} \di \dj \ln P
      &= \dci^n \di \dj \frac {\dci P}{P} \\
      &= x_i \dci^n \di \dj \frac {\dxi \dz P}{P} \\
      &= x_i \dci^n \di \dj \left(
           \dxi \frac {\dz P}{P} - \dz P · \dxi \frac 1 P
         \right) \\
      &= x_i \dci^n \di \dj \left(
           \dxi \frac {\dz P}{P} + \frac {\dz P}{P} · \frac {\dxi P}{P}
         \right) \\
      &= x_i \dxi \dz \underbrace {\dci^n \di \dj \ln P}_{=0} +
         x_i \dci^n \di \dj (\dz \ln P · \dxi \ln P) \\
      &= x_i \dci^n   (\di \dj \dz \ln P · \dxi \ln P +
                       \di \dz \ln P · \dj \dxi \ln P \\
      &\qquad \qquad + \dj \dz \ln P · \di \dxi \ln P +
                       \dz \ln P · \di \dj \dxi \ln P) \\
      &= 0
  \end {align*}
  (for the last step note that every summand has a factor that contains a
  $ \dt^2 \ln P $, $ \dz^2 \ln P $, or $ \dt \dz \ln P $ that gets at most $n$
  $ \dci $'s, so it vanishes by the induction assumption).

  \emph {Step 2.} By step 1 it suffices to prove the lemma in the case $ c=0 $.
  Note that then $P$ becomes a product of two terms of the form
    \[ R = \sum_k \frac {x^k}{k!} \; t^{ak} \qquad \mbox {and} \qquad
       S = \sum_k \frac {x^k}{k!} \; \prod_{i=0}^{bk-1} (z+t-i) \]
  where the first term contains all the $ x_i $ with $ (a_i,b_i)=(0,0) $ or $
  (a_i,b_i)=(1,0) $, and the second term all the $ x_i $ with $ (a_i,b_i)=(0,1)
  $. Obviously, it suffices to prove the lemma for $R$ and $S$ separately. But
    \[ R = \sum_k \prod_i \frac {(x_i t^a_i)^{k_i}}{k_i!}
         = \exp \Big( \sum_i x_i t^{a_i} \Big) \]
  and
    \[ S = \sum_k \frac {x^k}{k!} \; \binom {z+t}{\sum k} \; (\sum k) !
           = \Big( 1+\sum_i x_i \Big)^{z+t}, \]
  and in both cases it is obvious that the lemma holds.
\end {proof}

\begin {lemma} \label {tech-lemma-2}
  Let $ x_i $ be a collection of variables (possibly infinite), and let $ c_i
  \in \CC $. Define
    \[ Q(t) = \sum_k \frac {x^k}{k!} \; t \;
                \prod_{i=1}^{\sum k-1} (ck+t-i) \]
  in multi-index notation, where $k$ is a multi-index. Then $ \ln Q(t) $ is
  linear in $t$, i.e.\
    \[ (t\dt-1) \ln Q = 0. \]
\end {lemma}

\begin {proof}
  The proof is very similar to that of lemma \ref {tech-lemma}.

  \emph {Step 1.} We consider the $ c_i $ to be formal variables and show by
  induction on $n$ that for every $i$ and every $ n \ge 0 $
    \[ \mbox {if} \quad (t\dt-1) \ln Q \at = 0
       \quad \mbox {then} \quad \dci^n (t\dt-1) \ln Q \at = 0. \]
  So assume that $ \dci^j (t\dt-1) \ln Q \at = 0 $ for $ j \le n $. By
  definition of $Q$ we have $ \dci Q = x_i \dxi (\dt - \frac 1 t) Q $. Hence it
  follows that (everything in the following calculation is evaluated at $ c_i
  = 0 $):
  \begin {align*}
    \dci^{n+1} (t\dt-1) \ln Q
      &= \dci^n (t\dt-1) \frac {x_i \dxi (\dt-\frac 1 t) Q}{Q} \\
      &= x_i (t\dt-1) \left(
           \underbrace {\dci^n \left( \dt-\frac 1 t \right) \dxi \ln Q}_{=0} +
           \dci^n \dt \ln Q · \dxi \ln Q
         \right) \\
      &= x_i \dci^n \left(
           \dt \ln Q · \dxi (t\dt-1) \ln Q + \dt (t\dt-1) \ln Q · \dxi \ln Q
         \right) \\
      &= 0
  \end {align*}
  (for the last step note that every summand has a factor that contains a
  $ (t\dt-1) \ln Q $ that gets at most $n$ $ \dci $'s, so it vanishes by the
  induction assumption).

  \emph {Step 2.} By step 1 it suffices to prove the lemma in the case $ c=0 $.
  But then
    \[ Q(t) = \sum_k \frac {x^k}{k!} \; \prod_{i=0}^{\sum k-1} (t-i)
            = \Big( 1+\sum_i x_i \Big)^t, \]
  which obviously satisfies the statement of the lemma.
\end {proof}

  \end {appendix}

  % Main file: latex main


\begin {thebibliography}{XXXX}

\bibitem [B]{B} A. Bertram, \emph {Another way to enumerate rational curves
  with torus actions}, \preprint {math.AG}{9905159}.

\bibitem [CKYZ]{CKYZ} T. Chiang, A. Klemm, S. Yau, E. Zaslow, \emph {Local
  mirror symmetry: calculations and interpretations}, \preprint
  {hep-th}{9903053}.

\bibitem [Ga]{Ga} A. Gathmann, \emph {Absolute and relative Gromov-Witten
  invariants of very ample hypersurfaces}, \preprint {math.AG}{9908054}.

\bibitem [Ge]{Ge} E. Getzler, \emph {Topological recursion relations in genus
  2}, in \emph {Integrable systems and algebraic geometry} (Kobe/Kyoto,
  1997), 73--106, World Sci.\ Publishing, River Edge, NJ, 1998.

\bibitem [Gi]{Gi} A. Givental, \emph {Equivariant Gromov-Witten invariants},
  Internat.\ Math.\ Res.\ Notices \textbf {13} (1996), 613--663.

\bibitem [K]{K} B. Kim, \emph {Quantum hyperplane section principle for
  concavex decomposable vector bundles}, J. Korean Math.\ Soc.\ \textbf {37}
  (2000), no.\ 3, 455--461.

\bibitem [L]{L} Y. Lee, \emph {Quantum Lefschetz hyperplane theorem}, \preprint
  {math.AG}{0003128}.

\bibitem [LLY]{LLY} B. Lian, K. Liu, S. Yau, \emph {Mirror principle I}, Asian
  J. of Math.\ \textbf {1} (1997), no.\ 4, 729--763.

\bibitem [P]{P} R. Pandharipande, \emph {Rational curves on hypersurfaces
  (after A. Givental)}, Astérisque \textbf {252} (1998), Exp.\ No.\ 848, 5,
  307--340.

\bibitem [T]{T} N. Takahashi, \emph {Log mirror symmetry and local mirror
  symmetry}, \preprint {math.AG}{0004179}.

\end {thebibliography}
\end {document}